%% file: rkf.tex
\documentclass [a4paper, 12pt] {amsart}

\include {lwen}

\begin {document}

\title {Remarks on the Riesz-Kantorovich formula}
\author {D.~V.~Rutsky}
\email {rutsky@pdmi.ras.ru}
\date {\today}
\address {St.Petersburg Department
of Steklov Mathematical Institute RAS
27, Fontanka
191023 St.Petersburg
Russia}

\keywords {Riesz-Kantorovich formula, Riesz Decomposition Property, Hahn-Banach Theorem, M. Riesz Extension Theorem}

\begin {abstract}
The Riesz-Kantorovich formula expresses (under certain assumptions)
the supremum of two operators
$S, T : X \to Y$ where $X$ and $Y$ are ordered linear spaces as
$$
S \vee T (x) = \sup_{0 \leqslant y \leqslant x} [S (y) + T (x - y)].
$$
We explore some conditions that are sufficient for its validity, which enables us to get extensions of known results
characterizing lattice and decomposition properties of certain spaces of order bounded linear operators between $X$ and $Y$.
\end {abstract}

\maketitle

\setcounter {section} {-1}

\section {Introduction}

\label {anintroduction}


Ordered linear spaces, and especially Banach lattices, and various classes of linear operators between them,
play an important part in modern functional analysis (see, e.~g., \cite {positiveoperators})
and have significant applications to mathematical economics (see, e. g., \cite {aliprantisetal2000ec}, \cite {aliprantisetal2000}).
The following famous theorem was established independently
by M. Riesz in \cite {riesz1940}
and L.~V. Kantorovich
in \cite {kantorovich1936}
(all definitions are given in Section~\ref {ovl} below; we quote the statement from \cite {aliprantistourky2002}).
\begin {theoremrieszkantorovich}
If $L$ is an ordered vector space with the Riesz Decomposition Property and $M$ is a Dedekind complete vector lattice,
then the ordered vector space $\blinop {L} {M}$ of order bounded linear operators between $L$ and $M$
is a Dedekind complete vector lattice.  Moreover, if $S, T \in \blinop {L} {M}$, then
for each $x \in L_+$ we have
\begin {gather}
\label {stvee}
[S \vee T] (x) = \sup \{S y + T (x - y) \mid 0 \leqslant y \leqslant x \},
\\
\label {stwedge}
[S \wedge T] (x) = \inf \{S y + T (x - y) \mid 0 \leqslant y \leqslant x \}.
\end {gather}
\end {theoremrieszkantorovich}
In particular, this theorem implies that under its conditions every order bounded operator $T : L \to M$ is regular.
Equations \eqref {stvee} and \eqref {stwedge} and their variants describing the lattice operations of $\blinop {L} {M}$
are known as the \emph {Riesz-Kantorovich formulae}.  Related to these formulae is the following long-standing problem
that was posed in a number of slightly different versions (see, e.~g., \cite [\S 2] {positiveoperators}, \cite {aliprantistourky2002}):
if a supremum of two linear opeartors $S, T \in \mathcal L$ exists in a linear space $\mathcal L$ of linear operators between $L$ and $M$
(where of special interest are the space of order bounded operators and the space of regular operators), is it given by a Riesz-Kantorovich formula?
Although no satisfactory and definite answer seems to have been given yet in general,
significant progress has been made over the years and a number of partial answers had been found; see, e.~g.,
\cite {abramovich1990}, \cite {abramovichgeyler1982}, \cite {abramovichwickstead1991}, \cite {abramovichwickstead1993},
\cite {abramovichwickstead1997}, \cite {rooij1985}, \cite {aliprantisetal2000}, \cite {aliprantistourky2002}.
In remarkable works \cite {aliprantisetal2000} and \cite {aliprantistourky2002} (see also the monograph \cite {conesandduality})
motivated by recent advances in
welfare economics,
this question was answered affirmatively in the case of linear continous functionals (i.~e. $M = \mathbb R$) under relatively mild restrictions:
the order intervals of $L$ are assumed to be weakly compact, and either $L_+$ has an interior point or $L$ is an ordered
Fréchet space.  A new theory of super order duals, later generalized to a theory of
$\mathcal K$-lattices, was developed for this purpose (see \cite [\S 5.1] {conesandduality}).
The idea in a nutshell is to construct suitable extensions of the order duals and study their properties.
The earlier work \cite {aliprantisetal2000}
treating only the case of $L_+$ having an internal point (or, equivalently, an order unit), however,
used a somewhat tricky but quite elementary approach.

In the present work we extend the results of \cite {aliprantisetal2000} and \cite {aliprantistourky2002}.
A direct approach, inspired by the well-known characterization of convex functions as supremums of sets of affine functions,
allows us to get the complete answer in the case of order bounded functionals.
Unfortunately, this method does not seem to give much information regarding
the general case of order bounded linear operators acting from an ordered vector space $L$ into a Dedekind complete Riesz space $M$,
which is geometrically far more challenging compared to the case of functionals $M = \mathbb R$.
However, the same pattern still works to an extent if we assume that $L_+$ has an internal point and engage a natural extension of the Hahn-Banach theorem
to construct suitable majorants.
Although it only shows in this case that the Riesz-Kantorovich transform is linear on the set of internal points of $L_+$,
leaving under question its linearity on the entire $L_+$, this is still sufficient for many applications.
Moreover, this approach also allows us to characterize (under some additional assumptions)
a generalization of the Riesz Decomposition Property of $L$ in terms of the Riesz Decomposition Property
of $\blinop {L} {M}$ (as it was done, e.~g., in \cite {ando1962} for the duals of certain ordered Banach spaces)
without much effort and without embedding the spaces into their second duals.

We also give a characterization of conditions under which the Riesz-Kantorovich formulae are linear in terms of a
generalization of the Riesz Decomposition Property and investigate certain conditions under which this generalization
coincides with the (classical) Riesz Decomposition Property.

The paper is organized as follows.
Section~\ref {ovl} contains all the basic definitions and facts that we will need in the present work,
including the natural extension of the Hahn-Banach theorem mentioned above.
In Section~\ref {lrdps} we introduce a generalization of the Riesz Decomposition Property
and disuss some of its properties and some cases when it is equivalent to the (usual) Riesz Decomposition Property.
In Section~\ref {trk} we give a concise exposition of the Riesz-Kantorovich transform and investigate how its linearity is related to the
generalization of the Riesz Decomposition Property discussed in the previous section.
Section~\ref {dofsm} contains some of the more technically challenging parts of the proofs of the main results
that are presented in Section~\ref {rkf}.

To conclude the introduction,
let us note that although the condition that $L_+$ has an interior point
under which we give some answers regarding the case of the
general space $M$ is a fairly strong assumption (e.~g. among the classical Banach lattices $\lclassg {p}$ only
$\lclassg {\infty}$ satisfies it for all measures),
it seems plausible nonetheless that at least under some suitable conditions it can be relaxed to the assumption
that $L_+$ has a quasi-interior point (or a weak order unit; see, e.~g., \cite {borweinlewis1992}, \cite [Theorem 4.85] {positiveoperators1985}),
which is often satisfied.




\section {Ordered linear spaces}

\label {ovl}

In this section we introduce the definitions and basic properties that are used in the present work.
For more detail see e.~g. \cite {conesandduality}.
An \emph {ordered linear space} $X$, also called an \emph {ordered vector space}, is a linear space
(we only consider spaces over the field of reals $\mathbb R$) equipped with a partial order
$\geqslant$ which is consistent with the linear structure, that is for all $x, y, z \in V$ such that $x \geqslant y$
and $\lambda \geqslant 0$ we have $x + z \geqslant y + z$ and $\lambda x \geqslant \lambda y$.
The order $\geqslant$ is completely determined by the \emph {nonnegative cone} $X_+ = \{x \in X \mid x \geqslant 0\}$
(which is a convex pointed cone) in the sense that every convex cone $K \subset X$ of a linear space $X$ such that $K \cap (-K) = \{0\}$
determines a unique order ($x \geqslant y$ if and only if $x - y \in K$) that turns $X$ into an ordered linear space satisfying $X_+ = K$.
The field of reals $\mathbb R$ is itself an ordered linear space with respect to the natural order
of~$\mathbb R$.  We say that a cone $X_+$ is \emph {generating in $X$} if $X = X_+ - X_+$, that is for any $x \in X$ there exist some
$y, z \in X_+$ such that $x = y - z$.

Let $X$ be an ordered linear space.
An \emph {order interval} between a couple of points $x, y \in X$ such that $y \geqslant x$ is the set
$[x, y] = \{ z \in X \mid x \leqslant z \leqslant y \}$.  We say that a set $A \subset X$ is \emph {order bounded} if
$A \subset [x, y]$ for some $x, y \in X$, $y \geqslant x$.  If $y \geqslant x$ for all $x \in A$ then $y$ is called an \emph {upper bound} of $A$;
we sometimes denote it by $y \geqslant A$.
If there exists some $x \in X$ such that $y \geqslant x$ for any upper bound $y$ of $A$, i.~e. if $x$ is the \emph {least upper bound} of $A$,
then $x$ is called the \emph {supremum} of $A$ and we denote this
fact by $x = \sup A$.  Likewise, $x$ is called an \emph {infinum} of $A$, $x = \inf A$, if $x$ is the greatest lower bound.
Naturally, it is equivalent to existence of the supremum $-x$ of the set $-A = \{-y \mid y \in A\}$, and we have $\inf A = -\sup (-A)$
in this case.

Let $X$ and $Y$ be linear spaces.  We denote by $\linop {X} {Y}$ the set of all linear operators $T : X \to Y$.
If, additionally,
$X$ and $Y$ are equipped with topologies turning them into topological linear spaces,
we denote by $\clinop {X} {Y} \subset \linop {X} {Y}$ the set of all continuous operators.
Now suppose that $X$ and $Y$ are ordered linear spaces.
We say that a (not necessarily linear) map
$T : X \to Y$ is \emph {order bounded} if its values $T (A) = \{ T x \mid x \in A \}$ on order bounded sets $A \subset X$ are
order bounded.
$T$ is called \emph {positive} if $T x \geqslant 0$ for all $x \in X_+$; we denote this fact by $T \geqslant 0$.
$T$ is called \emph {strictly positive} if $T x > 0$ for all $x \in X_+ \setminus \{0\}$.
We say that $S \geqslant T$ for some $S : X \to Y$ if $S - T \geqslant 0$; this defines a natural partial order
on the set of all maps $T : X \to Y$.
$T$ is called \emph {regular} if $T \in \linop {X} {Y}$ and
$T = U - V$ for some positive $U, V \in \linop {X} {Y}$.
We denote the set of order bounded linear operators between ordered linear spaces $X$ and $Y$ by $\blinop {X} {Y}$,
and likewise the set of regular operators between $X$ and $Y$ by $\rlinop {X} {Y}$.
It is customary to use the notation $\orderdual {X}$ for the \emph {order dual} $\blinop {X} {\mathbb R}$ of an ordered linear space $X$;
however, since the focus of the present work is spaces of general order bounded operators
$\blinop {X} {Y}$
this notation might be a bit confusing, so we will avoid it.
It is easy to see that every positive linear operator is order bounded
and thus $\rlinop {X} {Y} \subset \blinop {X} {Y}$; this inclusion may be proper in general.
Both of these sets of operators may be very narrow: for example, if $X = \mathbb R^2$ is the plane with the lexicographic
order $X_+ = \left\{ (x, y) \mid \text {either $y > 0$ or $y = 0$ and $x \geqslant 0$} \right\}$
then $\rlinop {X} {\mathbb R} = \blinop {X} {\mathbb R} = \{f \mid f (x, y) = \alpha y \}$, and in general
ordered linear spaces with the lexicographic order have one-dimensional duals $\rlinop {X} {\mathbb R} = \blinop {X} {\mathbb R}$ irrespective to
the (finite or infinite) dimension of $X$.

Suppose that $X$ and $Y$ are ordered linear spaces.  A map $T : X_+ \to Y$ is called \emph {superadditive}
if $T (x + y) \geqslant T (x) + T (y)$ for all $x, y \in X_+$; $T$ is called \emph {additive} if
$T (x + y) = T (x) + T (y)$ for all $x, y \in X_+$.  If, in addition, $T$ is \emph {positively homogeneous},
that is $T (\lambda x) = \lambda T (x)$ for all $\lambda \geqslant 0$ and $x \in X_+$, we say that $T$ is
\emph {superlinear} and \emph {linear} respectively.
If the cone $X_+$ is generating in $X$ then any linear map $T : X_+ \to Y$
can be uniquely extended to a linear map on all $X$ by setting $T x = T y - T z$ for $x = y - z$, $x \in X$, $y, z \in X_+$.
We denote the set of all superlinear maps $T : X_+ \to Y$ by \suplinop {X_+} {Y}.

An ordered linear space $X$ is said to satisfy the \emph {Riesz Decomposition Property} if $[0, x] + [0, y] = [0, x + y]$
for all $x, y \in X_+$;
here and elsewhere a sum of sets $A$ and $B$ in a linear space is the Minkowski sum $A + B = \{a + b \mid a \in A, \,\, b \in B\}$.
It is well known that the Riesz Decomposition Property is equivalent to the \emph {interpolation property}: for any finite $A, B \subset X$ such that
$A \leqslant B$ in the sense that $a \leqslant b$ for all $a \in A$ and $b \in B$ there exists some $x \in X$ such that
$A \leqslant x \leqslant B$, i.~e. $a \leqslant x \leqslant b$ for all $a \in A$ and $b \in B$.
An ordered linear space $X$ is said to be a \emph {lattice} or \emph {Riesz space} if for any $x, y \in X$ the supremum $x \vee y = \sup \{x, y\}$
exists in $X$.  If this is the case then infinums $x \wedge y = -[(-x) \vee (-y)]$, positive parts $x^+ = x \vee 0$
and moduli $|x| = x \vee (-x)$ also exist.
A lattice always has the Riesz Decomposition Property.
An ordered linear space $X$ is called \emph {Dedekind complete} (or just \emph {complete}) if any set $A \subset X$ that has an upper bound also
has a supremum 
$\sup A$.
We will need the following natural extension of a well-known fundamental result in analysis.
\begin {theoremhahnbanach}
Suppose that $X$ and $Y$ are ordered linear spaces,
$Y$ is Dedekind complete, $L \subset X$ is a linear space such that $X_+ + L = X$,
$p : X_+ \to Y$ is a superlinear map and $T_0 : L \to Y$ is a linear operator satisfying
$T_0 x \geqslant p (x)$ for all $x \in L \cap X_+$.  Then there exists an extension $T \in \linop {X} {Y}$ of $T_0$ satisfying
$T x \geqslant p (x)$ for all $x \in X_+$.
\end {theoremhahnbanach}
This theorem is nearly always formulated with sublinear majorants $q = -p$, but superlinear minorants
suit the purposes of the present work better.
Setting $X_+ = X$ (i.~e. considering the full order, that is $x \geqslant y$ for all $x, y \in X$;
the assumption that $X_+$ is a partial order is not actually used)
yields what is often called the
Hahn-Banach-Kantorovich theorem (\cite {kantorovich1935}); letting $p = 0$ produces the
M.~Riesz extension theorem (\cite {riesz1923}).
Let us give a sketch of the proof. 
By the standard inductive reasoning it is sufficient to consider the case
$\dim X \slash L = 1$.  We fix some $x_0 \in X \setminus L$;
then $X = \{x + s x_0 \mid x \in L, s \in \mathbb R\}$ and
we need to define a value $y_0 \in Y$ in such a way that an extension
$T (x + s x_0) = T_0 x + s y_0$ would satisfy the necessary estimate
$T (x + s x_0) \geqslant p (x + s x_0)$ whenever $x + s x_0 \in X_+$.
Considering different signs of $s$ we see that this is equivalent to
\begin {equation}
\label {hbk}
T_0 u - p (u - x_0) \geqslant y_0 \geqslant p (v + x_0) - T_0 v
\end {equation}
for all $u, v \in L$ satisfying $u - x_0 \in X_+$ and $v + x_0 \in X_+$.
Observe that such $u$ and $v$ always exist since
by the assumptions
$$
X_+ + L = X = -X = -(X_+ + L) = -X_+ - L = -X_+ + L,
$$
and for any such $u$ and $v$ we have $u + v = (u - x_0) + (v + x_0) \in X_+ \cap L$, so
$$
T_0 u + T_0 v = T_0 (u + v) \geqslant p (u + v) = p ([u - x_0] + [v + x_0]) \geqslant p (u - x_0) + p (v + x_0)
$$
by superlinearity of $p$, which implies that
\begin {equation}
\label {hbk1}
T_0 u - p (u - x_0) \geqslant p (v + x_0) - T_0 v
\end {equation}
for any $u, v \in L$ satisfying $u - x_0 \in X_+$ and $v + x_0 \in X_+$.
From \eqref {hbk1} it is clear that the set of values of the right-hand part of~\eqref {hbk}
has an upper bound, and the least upper bound of the right-hand part is also a lower bound for the left-hand part.
Thus we may set
$$
y_0 = \sup \{p (v + x_0) - T_0 v \mid v + x_0 \in X_+\}
$$
and conclude that $y_0$ satisfies the required estimates~\eqref {hbk}.

Let $X$ be a linear space.
We say that $e \in A$ is an \emph {internal point} of a set $A \subset X$
(sometimes also called a \emph {core point} or an \emph {algebraically interior point})
if for any $x \in A$ there exists some $\varepsilon > 0$ such that $e + \alpha x \in A$ for all values $\alpha \in \mathbb R$
satisfying $-\varepsilon < \alpha < \varepsilon$.  We denote the set of internal points of $A$ by $\tint A$.
It is well known that this notion of an internal point has a lot in common with the usual topological notion of an interior point.
With only a limited number of caveats the \emph {line-open sets} $A$,
i.~e. the sets $A \subset X$ such that $A = \tint A$, behave very much like open sets in a Hausdorff linear topological
space.
A convex combination $(1 - \theta) y + \theta x$, $0 < \theta < 1$, of an internal point $x \in \tint A$ of $A$ and some $y \in A$ is also an internal point
of $A$.

Now let $X$ be an ordered linear space.
A vector $e \in X_+$ is called an \emph {order unit} of $X$ if for any $x \in X$ there exists
some $\lambda > 0$ such that $x \leqslant \lambda e$.  It is well known that $e$ is an order unit of $X$ if and only if $e$ is an internal point of
$X_+$.
\begin {comment}
As we will see later in Section~\ref {rkf}, Theorem~\ref {onedimext} implies positive answer to the
main question of this paper.
A more immediate consequence (obtained by setting $p = 0$)
is that for any order unit $e$ of $X$ there exists a positive functional
$f$ satisfying $f (e) = 1$.  This leads to the following result concerning positive functionals.
\begin {theorem}
\label {posfuncded}
Suppose that $X$ is a Dedekind complete lattice and $X_+$ is a generating cone in $X$.  Then for any $x \in X_+$ there exists a positive
linear functional $f : X \to \mathbb R$ satisfying $f (x) = 1$.
\end {theorem}
Indeed, let $L_x = \bigcup_{\lambda > 0} \lambda [-x, x]$ be the ideal in $X$ generated by $x$.
Then $x$ is an order unit of $L_x$ and
it follows from Theorem~\ref {onedimext} that there
exists a positive linear functional $f_0 : L_x \to \mathbb R$ satisfying $f_0 (x) = 1$.  We need to extend $f_0$ from $L_x$ to the entire $X$.
This can be done by means of order projection: let $f (z)  = \sup \{ f_0 (y) \mid y \in L_x \cap [0, z] \}$ for $z \in X_+$.
Then $f (z) \geqslant 0$ for all $z \in X_+$ and $f (z) = f_0 (z)$ for all $z \in L_x \cap X_+$ by positivity of $f_0$.
Observe that $X$ is a lattice and thus has the Riesz Decomposition Property, so for any $a, b \in X_+$ we have $[0, a + b] = [0, a] + [0, b]$
and therefore
$$
L_x \cap [0, a + b] = L_x \cap \left([0, a] + [0, b]\right) \supset L_x \cap [0, a] + L_x \cap [0, b].
$$
On the other hand, if $y \in L_x \cap [0, a + b]$ then there exist $y_0 \in [0, a]$ and $y_1 \in [0, b]$ such that $y = y_0 + y_1$;
since $L_x$ is an ideal, we have $y_0, y_1 \in L_x$ and therefore $y \in L_x \cap [0, a] + L_x \cap [0, b]$.
So
$$
L_x \cap [0, a + b] = L_x \cap [0, a] + L_x \cap [0, b]
$$
for any $a, b \in L_x$; thus it is easy to see that the map
$T : z \mapsto L_x \cap [0, z]$ taking values in sets of $L_x$ bounded from above is linear on $X_+$ and conclude that
$z \mapsto f (z) = \sup T (z)$
is linear and therefore can be uniquely extended to a positive linear functional $f : X \to \mathbb R$ since the cone $X_+$ is generating in $X$.
The proof of Theorem~\ref {posfuncded} is complete.
\end {comment}

Let $X$ be an ordered linear space.
The \emph {order topology} of $X$ is the finest locally convex topology
that has a base of convex circled open neighbourhoods at $0$ that absorb every order interval of $X$;
this topology is Hausdorff if and only if $\blinop {X} {\mathbb R}$ separates points of $X$
(i.~e. if and only if $f (x) = 0$ for all $f \in \blinop {X} {\mathbb R}$ implies $x = 0$).
A detailed description of order topologies
can be found in \cite {namioka1957}, \cite {schaefer1958}; see also \cite [Section~6] {schaefertvs} or \cite [Section~2.8] {conesandduality}.
An important property of the order topology is that
$\blinop {X} {\mathbb R} = \clinop {X} {\mathbb R}$ (see, e.~g., \cite [Theorem~2.69] {conesandduality}).

\begin {comment}
One of the useful properties of positive maps is their \emph {monotoneity}, that is
if $g$ is a monotone map then $x \geqslant y$ implies $g (x) \geqslant g (y)$.
The following simple property for the case $Y = \mathbb R$ is repeatedly used throughout the paper.
\begin {proposition}
\label {monsuprel}
Suppose that $X$ and $Y$ are ordered linear spaces and $Y$ is Dedekind complete.
Let $g : X \to Y$ be a monotone map, and suppose that as set $A \subset X$ has a supremum $\sup A$ in $X$.
Then its image $g (A) = \{g (x) \mid x \in A\}$ under the action of $g$ also has a supremum $\sup g (A)$ and
we have the relation
\begin {equation}
\label {monsuprelf}
\sup g (A) \leqslant g (\sup A).
\end {equation}
\end {proposition}
Indeed, monotoneity of $g$ implies that for every $y \geqslant A$ we also have $g (y) \geqslant g (A)$;
in particular, $g (A) \leqslant g (\sup A)$, and therefore $g (A)$ is bounded from above,
$\sup g (A)$ exists in $Y$ and \eqref {monsuprelf} follows.
Simple examples show that the inequality in \eqref {monsuprelf} may be proper, and indeed it often is.
\end {comment}

\section {$\mathcal L$-Riesz Decomposition Property}

\label {lrdps}

In this section we introduce and discuss a generalization of the Riesz Decomposition Property
that will be used in Section~\ref {trk} below to characterize conditions under which the Riesz-Kantorovich transform of
suitable linear maps is linear.

Suppose that $X$ is a linear space, $Y$ is an ordered linear space, $A \subset X$ is a convex set, $x \in X$ and $f \in \linop {X} {Y}$.
We say that \emph {$f$ separates $x$ from $A$} if either $f (x) \leqslant f (A)$ or $f (x) \geqslant f (A)$, where $f (A) = \{f (y) \mid y \in A\}$.
Suppose also that $Y$ is Dedekind complete.
We say that \emph {$f$ strictly separates $x$ from $A$} if $f$ separates $x$ from $A$ and
the appropriate condition is strengthened as either $f (x) < \inf f (A)$ or $f (x) > \sup f (A)$.  In the case $Y = \mathbb R$
these properties coincide with the usual separation (either strict or non-strict) by a linear functional.

\begin {definition}
\label {wrdp}
Suppose that $X$ and $Y$ are ordered linear spaces, $Y$ is Dedekind complete and $\mathcal L$ is a linear subspace of $\linop {X} {Y}$.
We say that $X$ satisfies the $\mathcal L$-Riesz Decomposition Property if for any $x, y \in X_+$
and $z \in [0, x + y] \setminus ([0, x] + [0, y])$ point $z$ cannot be strictly separated from
$[0, x] + [0, y]$ by a map from $\mathcal L$.
\end {definition}

\begin {figure} [tbph]
\includegraphics [width=13cm] {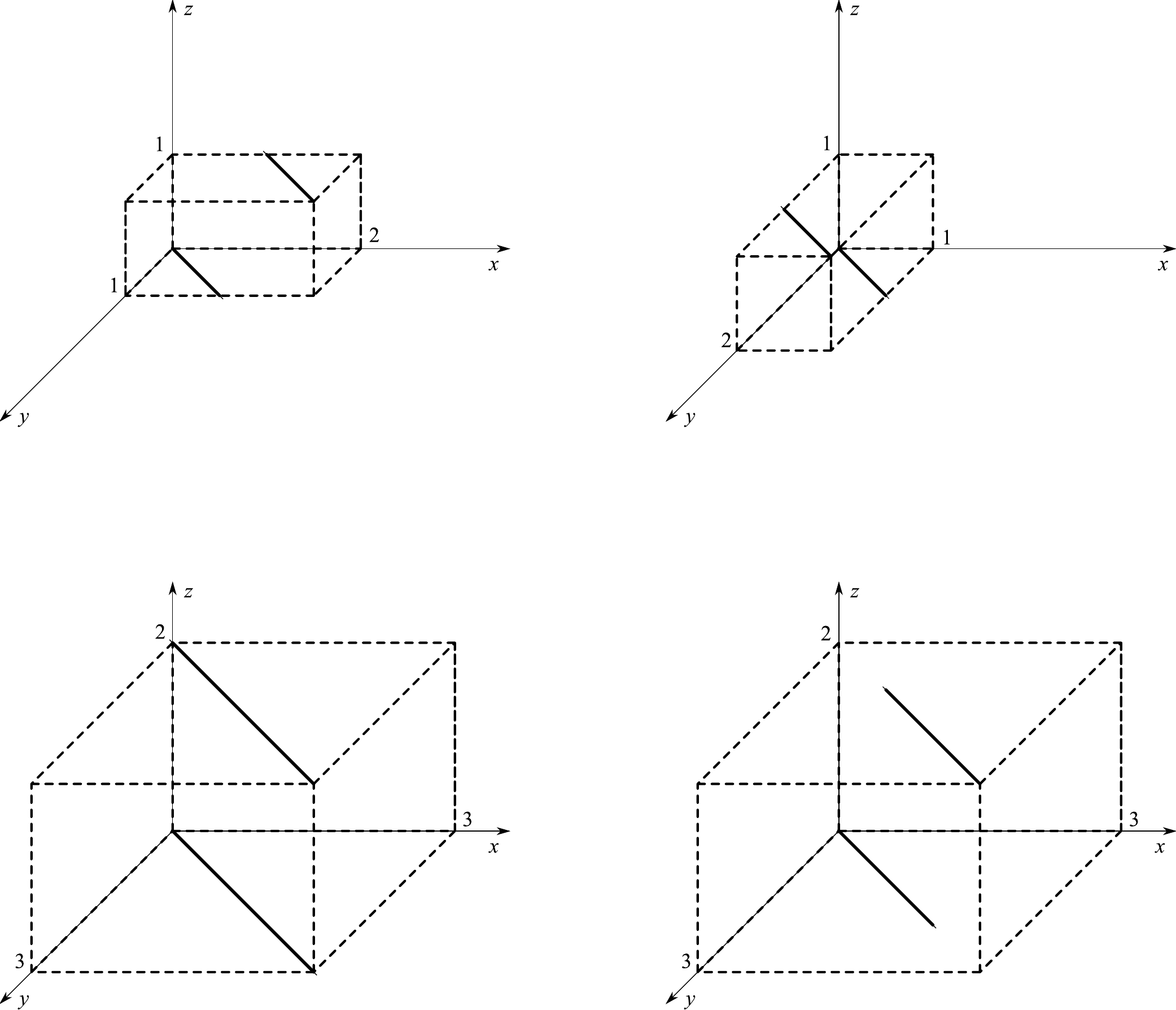}
\caption {Order intervals in an example of an order in $\mathbb R^3$ that satisfies the $\mathcal L$-Riesz Decomposition Property but does not satisfy the
Riesz Decomposition Property.
The order intervals $[0, (2, 1, 1)]$ and $[0, (1, 2, 1)]$ (top row) sum up to $[0, (2, 1, 1)] + [0, (1, 2, 1)]$ (bottom right); however,
the order interval $[0, (2, 1, 1) + (1, 2, 1)] = [0, (3, 3, 2)]$ (bottom left) has longer segments on the boundary.
}
\label {figrdpex}
\end {figure}

Trivially, if $X$ satisfies the Riesz Decomposition Property then it satisfies the $\mathcal L$-Riesz Decomposition Property with any
$\mathcal L$.
On the other hand, however, it is easy to construct an example that shows that the $\mathcal L$-Riesz Decomposition Property
is weaker than
the Riesz Decomposition Property for $Y = \mathbb R$ and the widest possible space $\mathcal L = \linop {X} {\mathbb R}$.
We will now give a detailed description of such an example; see Figure~\ref {figrdpex}.
Let
$X = \mathbb R^3$
and
$X_+ = \{ (x, y, z) \mid \text {$x, y, z > 0$ or $z = 0$, $x = y \geqslant 0$}\}$.
Observe that order interval $[0, (2, 1, 1)]$ is the union of an open rectangular box
$
\{(x, y, z) \mid 0 < x < 2, \,\, 0 < y < 1, \,\, 0 < z < 1\}
$
and two parallel semiopen segments connecting $(0, 0, 0)$ with $(1, 1, 0)$ not including $(1, 1, 0)$ 
and $(2, 1, 1)$ with $(1, 0, 1)$ not including $(1, 0, 1)$.  The order interval $[0, (1, 2, 1)]$ has similar structure because
it is the mirror image of $[0, (2, 1, 1)]$ with respect to the symmetry plane $\{(x, y, z) \mid x = y \}$ of $X_+$.  It is easy to see that
$
[0, (2, 1, 1) + (1, 2, 1)] = [0, (3, 3, 2)]
$
is the union of an open rectangular box
$$
\{(x, y, z) \mid 0 < x < 3, \,\, 0 < y < 3, \,\, 0 < z < 2\}
$$
and two parallel semiopen segments connecting $(0, 0, 0)$ with $(3, 3, 0)$ not including $(3, 3, 0)$ and
$(3, 3, 2)$ with $(0, 0, 2)$ not including $(0, 0, 2)$.  Thus
$[0, (3, 3, 2)] \setminus \left([0, (2, 1, 1)] + [0, (1, 2, 1)]\right)$ is a union of a semiopen segment between $(0, 0, 3)$ and $(1, 1, 3)$
not including $(0, 0, 3)$
and a semiopen segment between $(2, 2, 0)$ and $(3, 3, 0)$ not including $(3, 3, 0)$, and so
$X$ does not satisfy the Riesz Decomposition Property.  However, observe that the closure of the order intervals $[0, a]$, $a \in X_+$
(in the standard topology of $\mathbb R^3$, which is the only topology there is that turns $\mathbb R^3$ into a Hausdorff linear topological space)
coincides with the order intervals of $\mathbb R^3$ ordered by the standard cone
$$
\mathbb R^3_+ = \{(x, y, z) \mid x \geqslant 0, \,\, y \geqslant 0, \,\, z \geqslant 0 \}
$$
with which $\mathbb R^3$
is a lattice and thus has
the Riesz Decomposition Property.  Therefore for all $a, b \in X_+$ the set $[0, a + b] \setminus \left([0, a] + [0, b]\right)$ lies in the closure
of $[0, a + b]$.
Since $X$ is finite dimensional,
every functional from $\mathcal L$
is continuous and therefore cannot strictly separate a convex set from a point in its closure.  Thus
$X$ satisfies the $\mathcal L$-Riesz Decomposition Property for any $\mathcal L$.
Of course, if the set $\mathcal L$ is very small with respect to $\linop {X} {Y}$ in a certain sense, for example if $X$ is infinite-dimensional
and all operators from $\mathcal L$ factor through a finite-dimensional space, then the $\mathcal L$-Riesz Decomposition Property
may be \emph {much} more general than the Riesz Decomposition Property; an obvious example is that \emph {any} couple of ordered linear spaces
$X$, $Y$ satisfies the $\mathcal L$-Riesz Decomposition Property if all operators from $\mathcal L$ factor through a one-dimensional space,
or just $\mathcal L = \{0\}$.
On the other hand, it is easy to enumerate a number of cases where
both properties are actually equivalent based on various separation theorems; let us give a couple of examples.

\begin {proposition}
\label {rdplrdp}
Suppose that ordered linear space $X$ is equipped with a topology that turns it into a locally convex Hausdorff linear space and
order intervals $[0, x]$ are compact for any $x \in X_+$.
Then $X$ satisfies the Riesz Decomposition Property if and only if $X$ satisfies the $\clinop {X} {\mathbb R}$-Riesz Decomposition Property.
\end {proposition}
Indeed, suppose that $X$ satisfies the $\clinop {X} {\mathbb R}$-Riesz Decomposition Property under the conditions of Proposition~\ref {rdplrdp}
and there exist some $x, y \in X_+$ and $z \in [0, x + y] \setminus \left( [0, x] + [0, y] \right)$.
$[0, x] + [0, y]$ is a compact convex set because it is the Minkowski sum of two compact convex sets.  Then by a standard separation theorem
$z$ can be strictly separated from $[0, x] + [0, y]$ by some $f \in \clinop {X} {\mathbb R}$, a contradiction.

The next example concerns ordered linear spaces $X$ such that every point of $X_+$ is internal.
It is a bit more involved, and we seem to need additional assumptions.  First we need a suitable separation theorem.
Suppose that $X$ is an ordered linear space and $B \subset X$ is a convex set.
We say that a cone $X_+$ \emph {has a base $B$} if for any $x \in X_+ \setminus \{0\}$
there exist unique $y \in B$ and $\lambda > 0$ such that $x = \lambda y$;
for details on cones with bases see, e.~g., \cite [Section~1.7] {conesandduality}.
Existence of a base for $X_+$ is equiuvalent to existence of a strictly positive functional $f$ in $X$; in this case the set
$B = \{ x \in X_+ \mid f (x) = 1\}$ is a base for $X_+$.
Restrictions of maps from $\linop {X} {Y}$ on $B$ are exactly the affine functionals on $B$ taking values in $Y$, and restrictions
$K \cap B$ of convex cones $K \subset X_+$ are exactly the convex subsets of $B$.
\begin {definition}
\label {finap}
Let $X$ be an ordered linear space such that the cone $X_+$ has a base $B$.
We say that $X_+$ has bounded aperture on parallel lines if for any $x \in X$ there exists some $c_x \geqslant 0$ such that
for any $y \in B \setminus \{x\}$ the length of the intersection $B \cap \ell$ is at most $c_x$,
where $\ell$ is the line passing through points $y$ and $x + y$.
\end {definition}
It is sufficient to consider only points $x$ from the affine hull of $B$ in Definition~\ref {finap},
since otherwise $B \cap \ell$ is a single point.
Also, since all bases of a cone are affinely isomorphic
and affine maps preserve relative lengths of parallel segments, Definition~\ref {finap} does not depend on a particular choice of a base $B$.
It is easy to see that for the finite dimensional ordered spaces the cone $X_+$ has bounded aperture if and only if its closure
does not contain any lines, i.~e. one-dimensional spaces.
It is also easy to see that if $X$ is a Hausdorff ordered linear topological space such that the cone $X_+$ has a compact base $B$ then
the cone $X_+$ has bounded aperture on parallel lines; however, this property seems to be much weaker than any kind of compactness.
For example, if $B$ is a bounded convex set in an infinite dimensional
normed linear space then the cone $X_+$ has bounded aperture on parallel lines (since the length of $B \cap \ell$ in Definition~\ref {finap}
is proportional to $\|z_0 - z_1\|$, where $z_0$ and $z_1$ are the endpoints of $B \cap \ell$),
but at the same time $B$ may lack any meaningful kind of compactness.  However, it allows us to establish the following separation theorem.
\begin {proposition}
\label {regfuncsep}
Suppose that $X$ is an ordered linear space such that the cone $X_+$ is generating for $X$ and $X_+$ has bounded aperture on parallel lines.
Let $K_0, K_1 \subset X_+$ be two nonempty convex cones such that $K_0 \cap K_1 = \emptyset$ and
$\tint K_0 \neq \emptyset$.
Then there exists a regular functional $f \in \rlinop {X} {\mathbb R}$, $f \neq 0$, such that $f$ separates $K_0$ and $K_1$.
\end {proposition}
Indeed, let $B$ be a base of the cone $X_+$, let $L$ be the affine hull of $B$ and let $A_j = K_j \cap B$ for $j \in \{0, 1\}$.
Then $\tint A_0 \neq \emptyset$ as a subset of $L$ (for the rest of the proof we denote by $\tint D$ the sets of internal points of
$D \subset L$ in $L$ rather than in $X$),
and by the standard separation theorem there exists a linear functional
$h_{00} \in \linop {L} {\mathbb R}$ separating $A_0$ and $A_1$.  This means that for some $\alpha \in \mathbb R$ the hyperplane
$H = \{y \in L \mid h_{00} (y) = \alpha\}$ separates $A_0$ and $A_1$ in $L$.
Let $x_0 \in \tint A_0$ and let $y_0 \in H \cap B$; then $x_0 \notin H$.  Space $L$ can then be parametrized as
$L = \{s (x_0 - y_0) + y \mid s \in \mathbb R, y \in H \}$.
By the bounded aperture property assumption there exists some $c \in \mathbb R$ such that
\begin {equation}
\label {positiveonb}
B \subset \{s (x_0 - y_0) + y \mid |s| \leqslant c, \,\, y \in H \}.
\end {equation}
Let us define two affine functionals for $z = s (x_0 - y_0) + y$, $s \in \mathbb R$, $y \in H$ by
$h_0 (z) = 1 + \frac 1 {2 c} s$ and $h_1 (z) = 1 - \frac 1 {2 c} s$, and let $f_0, f_1 \in \linop {X} {\mathbb R}$
be the unique linear extensions of $h_0$ and $h_1$.
\eqref {positiveonb} implies that $h_0$ and $h_1$ are positive on $B$, so we have $f_0, f_1 \geqslant 0$ and thus $f = f_0 - f_1$
is a regular functional that is a unique extension of $h = h_0 - h_1$.  By construction $h$ is an affine functional that defines the same hyperplane
$H$, i.~e. $H = \{y \in L \mid h (y) = 0 \}$, and it is easy to see that we also have $h (A_0) \leqslant 0$ and $h (A_1) \geqslant 0$,
so $h$ separates $A_0$ and $A_1$, and therefore $f$ separates $K_0$ and $K_1$.  The proof of Proposition~\ref {regfuncsep} is complete.

It is not difficult to see that the assumption that $X_+$ has bounded aperture on parallel lines cannot be dropped
from Proposition~\ref {regfuncsep}: if the affine hull of the base $B$ of $X_+$ contains a line $\ell$,
then any positive functional $f_0$ has to be constant on $\ell$, and thus no
$f \in \rlinop {X} {\mathbb R}$, $f \neq 0$, can separate any cones $K_0, K_1 \subset X_+$ such that $\ell \cap \tint K_j \neq \emptyset$ for $j \in \{0, 1\}$.

We will need the following simple corollary of Proposition~\ref {regfuncsep}.
\begin {proposition}
\label {regfunchp}
Suppose that $X$ is an ordered linear space such that the cone $X_+$ is generating in $X$, $X_+$ has bounded aperture on parallel lines
and $\tint X_+ \neq \emptyset$.
Let $H \subset X$ be a hyperplane.
Then there exists a regular functional $f \in \rlinop {X} {\mathbb R}$ such that $H = \{x \mid f (x) = \alpha\}$ with some
$\alpha \in \mathbb R$.
\end {proposition}
First, fix any $x_0 \in H$, and let $H_0 = H - x_0$; we want to find some
$f \in \rlinop {X} {\mathbb R}$ such that $H_0 = \{x \mid f (x) = 0\}$.
Let $h_0 \in \linop {X} {\mathbb R}$ be a (not necessarily regular) linear functional defining $H_0$,
that is $H_0 = \{x \mid h_0 (x) = 0\}$.
If $H_0 \cap X_+ = \{0\}$ then either $h_0$ or $-h_0$ is positive and we can put $f = h_0$ right away.
Otherwise let
$$
K_0 = \{x \in X_+ \mid h_0 (x) < 0\}
$$
and
$$
K_1 = \{x \in X_+ \mid h_0 (x) > 0\}.
$$
By assumptions there exists some $z \in \tint X_+$.  We may assume that $z \notin H_0$, since otherwise we may replace $z$ with a convex combination
of $z$ and some point in $X_+ \setminus H_0$ (we cannot have $X_+ \subset H_0$ since $X_+$ is assumed to be generating in $X$).
Therefore
$$
z \in K_j = \{h_0 < 0\} \cap \tint X_+ = \tint K_j
$$
for some $j \in \{0, 1\}$;
by changing the sign of $h_0$ we may thus assume that $\tint K_0 \neq \emptyset$.
Applying Proposition~\ref {regfuncsep} yields a functional $f_0 \in \rlinop {X} {\mathbb R}$ separating $K_0$ and $K_1$,
and since the cone $X_+$ is generating in $X$, it is easy to see that $H_0 = \{ x \mid f (x) = 0\}$.
Finally, $H = \{x \mid f (x) = f (x_0)\}$, and the proof of Proposition~\ref {regfunchp} is complete.

\begin {proposition}
\label {aocsd}
Suppose that $X$ is a linear space and $B \subset A \subset X$ are nonempty, convex and line-open, i.~e. $A = \tint A$ and $B = \tint B$.
Then $A \setminus B \neq \emptyset$ implies that $\tint (A \setminus B) \neq \emptyset$;
moreover,
\begin {equation}
\label {aseplf}
A \cap \{z \mid g (z) > \alpha\} \subset \tint (A \setminus B)
\end {equation}
for some $g \in \linop {X} {\mathbb R}$ and $\alpha \in \mathbb R$.
\end {proposition}
Indeed, let $y_0 \in A \setminus B$ and $x_0 \in \tint B$.  By the standard separation theorem there exists a hyperplane
$H = \{ y \mid g (y) = \alpha\}$ separating $y_0$ and $B$.  By changing the sign of $g$ if necessary we may assume that $g (x_0) < \alpha$ and
$g (y_0) \geqslant \alpha$.
Since $A$ is line-open, there exists some $\lambda > 1$
such that $y = x_0 + \lambda (y_0 - x_0) \in A$.  Therefore $g (y) = g (y_0) + (\lambda  - 1) (g (y_0) - g (x_0))  > \alpha$.
This means that $y$ belongs to the line-open set \eqref {aseplf}.

We are now ready to give an equivalence result for ordered linear spaces with line-open cones.
\begin {theorem}
\label {openlrdp}
Suppose that $X$ is an ordered linear space satisfying $X_+ = \{0\} \cup \tint X_+$, the cone $X_+$ is generating for $X$,
$X_+$ has bounded aperture
on parallel lines
and $Y$ is an ordered linear space such that the cone $Y_+$ is generating for $Y$.
Let $\mathcal L \subset \linop {X} {Y}$ be a linear space of linear maps such that $y_0 \cdot \rlinop {X} {\mathbb R} \subset \mathcal L$
for some $y_0 \in Y$.
Then the following conditions are equivalent.
\begin {enumerate}
\item
$X$ has the Riesz Decomposition Property.
\item
$X$ has the $\mathcal L$-Riesz Decomposition Property.
\end {enumerate}
\end {theorem}
Indeed, $1 \Rightarrow 2$ is trivial and so we only need to establish $2 \Rightarrow 1$.
It is easy to see that we may, without loss of generality, assume that $Y = \mathbb R$ and $\rlinop {X} {\mathbb R} \subset \mathcal L$
under the conditions of Proposition~\ref {openlrdp}.
Suppose that $X$ has the $\mathcal L$-Riesz Decomposition Property but there exist some $x, y, z \in X_+$
such that
$z \in [0, x + y] \setminus \left([0, x] + [0, y]\right)$.
Since by the assumptions $X_+ = \{0\} \cup \tint X_+$, we have $[0, x] = \{0\} \cup \tint [0, x]$ and thus
$[0, x] + [0, y] = \{0\} \cup \tint \left( [0, x] + [0, y]\right)$.
It is easy to see that $\tint [0, x] \neq \emptyset$ unless $X$ is trivial.
Let $x_0 \in [0, x] + [0, y]$.
Using Proposition~\ref {aocsd}, we may assume that
\begin {equation}
\label {ztint}
z \in A = \tint \left( [0, x + y] \setminus \left([0, x] + [0, y]\right) \right) \subset [0, x + y] \cap \{z \mid g (z) > \alpha\}
\end {equation}
with some linear functional $g$.
Let $H = \{z \mid g (z) = \alpha\}$.
Then by Proposition~\ref {regfunchp} there exists a regular functional $f \in \rlinop {X} {\mathbb R}$
such that
$
H = \{a \mid f (a) = \beta\}
$
with some $\beta \in \mathbb R$.
This means that $f$ separates $A$ and $[0, x] + [0, y]$, and therefore strictly separates
$z \in A = \tint A$ from $[0, x] + [0, y]$, which contradicts the assumption that
$X$ satisfies the $\mathcal L$-Riesz Decomposition Property.  The proof of Theorem~\ref {openlrdp} is complete.

\section {The Riesz-Kantorovich transform}

\label {trk}

In this section we introduce the Riesz-Kantorovich transform and discuss some of its fundamental properties that will be needed
for the results in the following sections.
The arguments are mostly based on \cite {conesandduality}; we omit direct references.

\begin {definition}
\label {rktd}
Let $X$ and $Y$ be ordered linear spaces, $n \geqslant 1$ and $T = \left\{T_1, \ldots, T_n\right\}$ be a finite collection of superlinear maps
$T_j : X_+ \to Y$ for $1 \leqslant j \leqslant n$.
The Riesz-Kantorovich transform
$$
\rkts {T} = \rkts {T_1, \ldots, T_n} : x \mapsto \rkt {T} {x}
$$
is defined by
\begin {multline}
\label {rktdef}
\rkt {T} {x} = \rkt {T_1, \ldots, T_n} {x} =
\\
\sup \left\{ \sum_{j = 1}^n T_j (x_j) \mid x_j \in X_+ \text { for $1 \leqslant j \leqslant n$ and } \sum_{j = 1}^n x_j \leqslant x \right\}
\end {multline}
for $x \in X_+$ whenever this supremum exists in $Y$.
The positive Riesz-Kantorovich transform is defined by
$$
\rktp {T} {x} = \rktp {T_1, \ldots, T_n} {x} = \rktp {0, T_1, \ldots, T_n} {x}.
$$
\end {definition}
It is easy to see that if maps $T_j$ are order bounded and $Y$ is Dedekind complete then the Riesz-Kantorovich transform $\rkt {T} {x}$ is well-defined
and $\rkts {T_1, \ldots, T_n} \geqslant T_j$ for all $1 \leqslant j \leqslant n$.
Since
$[0, x] + [0, y] \subset [0, x + y]$ for all $x, y \in X_+$, it is easy to verify that $\rkts {T}$ is superadditive.
Since the Riesz-Kantorovich transform is also positively homogeneous,
we see that $\rkts {T}$ is a superlinear map that majorizes $T_j$ for all $1 \leqslant j \leqslant n$.  Observe also that
$\rkts {T}$ is order bounded whenever the operators from collection $T$ are.

A remarkable property of the Riesz-Kantorovich transform is that, whenever it is well defined,
it gives the infinum 
in the set of all superlinear maps.
\begin {proposition}
\label {rktinf}
Suppose that under the assumptions of Definition~\ref {rktd} map $\rkts {T}$ is well-defined.
Then
\begin {multline}
\label {rktinfrel}
\rkts {T} =
\\
\inf \{S : X_+ \to Y \mid \text {$S$ is superlinear and $S \geqslant T_j$ for all $1 \leqslant j \leqslant n$}\}.
\end {multline}
In particular, $\rkts {T} = \bigvee_{j = 1}^n T_j$ in \suplinop {X_+} {Y}.
\end {proposition}
Indeed, suppose that under the assumptions of Proposition~\ref {rktinf}
map $S : X_+ \to Y$ is superlinear and $S \geqslant T_j$ for all $1 \leqslant j \leqslant n$.
Then for any $x \in X_+$ and $x_j \in X_+$, $1 \leqslant j \leqslant n$, such that $\sum_{j = 1}^n x_j \leqslant x$, we have
\begin {equation}
\label {sssup}
S (x) \geqslant \sum_{j = 1}^n S (x_j) \geqslant \sum_{j = 1}^n T_j (x_j)
\end {equation}
by superadditivity of $S$, and \eqref {sssup} implies $S \geqslant \rkts {T}$.

As an immediate application of Proposition~\ref {rktinf}, let us show that the Riesz-Kantorovich transform is associative.
\begin {proposition}
\label {rkass}
Suppose that $X$ and $Y$ are ordered linear spaces and $Y$ is Dedekind complete.
Then for any order bounded superlinear maps $Q, R, S : X \to Y$ we have the relations
\begin {equation}
\label {asseq}
\rkts {\rkts {Q, R}, S} = \rkts {Q, R, S} = \rkts {Q, \rkts {R, S}}.
\end {equation}
\end {proposition}
By symmetry it is sufficient to establish the first equality in \eqref {asseq}.
Observe that \eqref {rktinfrel} implies that $\rkts {Q, R, S} \geqslant \rkts {Q, R}$, and another application of Proposition~\ref {rktinf}
yields $\rkts {Q, R, S} \geqslant \rkts {\rkts {Q, R}, S}$.
Let us verify the converse inequality.  Suppose that $w \in X_+$.  Then for any $x, y, z \in X_+$ such that $x + y + z \leqslant w$
we have
$$
\rkt {\rkts {Q, R}, S} {w} \geqslant \rkt {Q, R} {x + y} + S (z) \geqslant Q (x) + R (y) + S (z),
$$
so $\rkts {\rkts {Q, R}, S} \geqslant \rkts {Q, R, S}$.  It is easy to see that the Riesz-Kantorovich transform of any finite number of maps
is stable under arbitrary association of terms in the spirit of Proposition~\ref {rkass}.

Let us now briefly consider the question whether the inequality in the definition of \eqref {rktdef}
can be replaced by an equality.  A moment's reflection shows that it is generally not the case.
For example, let $X = \mathbb R^3$ with the standard lattice order
$X_+ = \left\{ (x, y, z) \mid x, y, z \geqslant 0 \right\}$ that appeared in Section~\ref {lrdps} above,
and let $f ((x, y, z)) = -x + y - z$ and $g ((x, y, z)) = x - y - z$.
Contemplating the extreme points of the order interval $[0, (x, y, z)]$ relative to the gradients
$(-1, 1, -1)$ and $(1, -1, 1)$ of $f$ and $g$, one immediately notices that
the supremum in \eqref {rktdef} for $\rkt {f, g} {(x, y, z)}$, $(x, y, z) \in X_+$, is attained
at a decomposition $(0, y, 0) + (x, 0, 0)$, which is strictly smaller than $(x, y, z)$ in the order of $X$ unless $z = 0$,
and that the corresponding supremum value
$$
\rkt {f, g} {(x, y, z)} = x + y
$$
would not have been attained had we replaced
inequality in the definition \eqref {rktdef} by equality.
However, it is easy to see that if one of the maps $T_j$ is positive then the supremum in \eqref {rktdef}
is attained at some decomposition $\sum_{j = 1}^n x_j = x$.
\begin {proposition}
\label {rkteq}
Let $X$ and $Y$ be ordered linear spaces, $n \geqslant 1$ and let $T = \left\{T_1, \ldots, T_n\right\}$ be a finite collection of linear maps
$T_j : X_+ \to Y$ for $1 \leqslant j \leqslant n$.
Suppose also that $T_l \geqslant 0$ for some $1 \leqslant l \leqslant n$.
Then
\begin {multline}
\label {rktdefeq}
\rkt {T} {x} =
\sup \left\{ \sum_{j = 1}^n T_j (x_j) \mid x_j \in X_+ \text { for $1 \leqslant j \leqslant n$ and } \sum_{j = 1}^n x_j = x \right\}.
\end {multline}
\end {proposition}
Indeed, the value of the right-hand part of \eqref {rktdefeq} is less or equal than $\rkt {T} {x}$ because the supremum is taken over a smaller
set of decompositions.
On the other hand, let $\sum_{j = 1}^n x_j \leqslant x$ be
any decomposition
in the definition \eqref {rktdef}.  Then we can construct another decomposition $x = \sum_{j = 1}^n \tilde x_j$ by taking
$\tilde x_j = x_j$ for $j \neq l$ and $\tilde x_l = x_l + \left[x - \sum_{j = 1}^n x_j\right]$ so that
$
\sum_{j = 1}^n T_j (\tilde x_j) = \sum_{j = 1}^n T_j (x_j) + T_l \left( x - \sum_{j = 1}^n x_j \right) \geqslant 
\sum_{j = 1}^n T_j (x_j).
$
It follows that 
\begin {multline*}
\rkt {T} {x} \leqslant
\sup \left\{ \sum_{j = 1}^n T_j (x_j) \mid x_j \in X_+ \text { for $1 \leqslant j \leqslant n$ and } \sum_{j = 1}^n x_j = x \right\},
\end {multline*}
which establishes \eqref {rktdefeq}.

Propositions~\ref {rkass} and \ref {rktdefeq} allow us to reduce many questions regarding the general Reisz-Kantorovich transform
to the case of $\rktps {T}$ of a single map $T$.  Observe that \eqref {rktdefeq} implies that
$$
\rkts {Q, R} = \rkts {0, Q - R} + R = \rktps {Q - R} + R
$$
for a positive linear operator $R \geqslant 0$ and a superlinerar map $Q$. Combining this with Proposition~\ref {rkass}
leads to the following useful observation.
\begin {proposition}
\label {onlysing}
Let $X$ and $Y$ be ordered linear spaces and $\mathcal L \subset \linop {X} {Y}$ a linear space.
Suppose that for any $T \in \mathcal L$ the Riesz-Kantorovich transform $\rktps {T}$ is well-defined and
$\rktps {T} \in \mathcal L$.  Then for any finite collection of operators $T = \{T_1, \ldots, T_n\}$,
$T_j \in \mathcal L$, $1 \leqslant j \leqslant n$, such that $T_l \geqslant 0$ for some $1 \leqslant l \leqslant n$,
we also have $\rkts {T} = \rkts {T_1, \ldots, T_n} \in \mathcal L$.
\end {proposition}

Since every linear map is superlinear, Proposition~\ref {rktinf} implies that if the Riesz-Kantorovich transform $\rkts {T}$ of a finite collection
of \emph {linear} operators is well-defined and linear then $\rkts {T} = \bigvee_{j = 1}^n T_j$ in \linop {X_+} {Y} as well as in
\suplinop {X_+} {Y}.
This leads to the Riesz-Kantorovich theorem mentioned in Section~\ref {anintroduction} as soon as one
notices that the Riesz Decomposition Property of $X$ implies linearity of $\rktps {T}$ for any $T \in \blinop {X} {Y}$.
More generally, it is not difficult to see that the $\mathcal L$-Riesz Decomposition Property
introduced in Section~\ref {lrdps} above is, in many cases, essentially
equivalent to linearity of the Riesz-Kantorovich transform.
In the rest of the current section we
explore the relationship between $\mathcal L$-Riesz Decomposition Property and linearity of the Riesz-Kantorovich transform.
\begin {proposition}
\label {lrdptorkf}
Suppose that $X$ and $Y$ are ordered linear spaces,
$Y$ is a Dedekind complete lattice, $\mathcal L$ is a linear subspace of $\linop {X} {Y}$,
$\rlinop {X} {Y} \subset \mathcal L$
and $X$ has the $\mathcal L$-Riesz Decomposition Property.
Then for any finite collection of operators $T = \{T_1, \ldots, T_n\}$
such that $T_j \in \mathcal L$ for all $1 \leqslant j \leqslant n$ and $T_l \geqslant 0$ for some $1 \leqslant l \leqslant n$,
we also have $\rkts {T} \in \mathcal L$.
\end {proposition}
Indeed, by Proposition~\ref {onlysing} it is sufficient to establish that for any $T \in \mathcal L$ the positive Riesz-Kantorovich
transform $\rktps {T}$ is linear; since $\rktps {T} \geqslant 0$, that would imply $\rktps {T} \in \mathcal L$ by the assumptions.
Since $\rktps {T}$ is superlinear, we have
$\rktp {T} {x + y} \geqslant \rktp {T} {x} + \rktp {T} {y}$ for all $x, y \in X_+$.
Suppose that this inequality is strict for some $x, y \in X_+$.  Then
$y_0 = \rktp {T} {x + y} - \rktp {T} {x} + \rktp {T} {y} > 0$, and from the definition of the supremum
it follows that
there exists some $z \in [0, x + y]$ such that
\begin {equation}
\label {tseplrdp}
T z \geqslant \frac 1 2 y_0 + \rktp {T} {x} + \rktp {T} {y} \geqslant
\frac 1 2 y_0 + T u + T v
\end {equation}
for all $u \in [0, x]$ and $v \in [0, y]$.
But \eqref {tseplrdp} means that $T$ strictly separates $z \in [0, x + y]$ from $[0, x] + [0, y]$,
which contradicts the assumption that $X$ satisfies the $\mathcal L$-Riesz Decomposition Property.
This concludes the proof of Proposition~\ref {lrdptorkf}.

The converse to Proposition~\ref {lrdptorkf} is established by essentially the same reasoning.
\begin {proposition}
\label {lrdp}
Suppose that $X$ and $Y$ are ordered linear spaces and $\mathcal L$ is a linear subspace of $\linop {X} {Y}$.
Suppose also that for any $T \in \mathcal L$ map
$\rktps {T}$ is well-defined and linear.  Then $X$ satisfies the $\mathcal L$-Riesz Decomposition Property.
\end {proposition}
Suppose that, on the contrary, $X$ does not satisfy the
$\mathcal L$-Riesz Decomposition Property.  This means that there exist
some $x, y \in X_+$, $z \in [0, x + y] \setminus \left( [0, x] + [0, y]\right)$
and $T \in \mathcal L$
such that
$T$ strictly separates $z$ from $[0, x] + [0, y]$.
By changing the sign of $T$ if necessary we may assume that
\begin {equation}
\label {tztuv}
T z > \sup \{T (u + v) \mid u \in [0, x], v \in [0, y]\}.
\end {equation}
Let $y_0 = T z - \sup \{T (u + v) \mid u \in [0, x], v \in [0, y]\} > 0$.
Then \eqref {tztuv} implies that
\begin {equation}
\label {tztuv2}
T z \geqslant y_0 + T u + T v
\end {equation}
for all $u \in [0, x]$ and $v \in [0, y]$.
Taking supremums of the right-hand part of \eqref {tztuv2} one after the other yields
\begin {multline*}
T z \geqslant y_0 + \sup \{T u \mid u \in [0, x]\} +
\sup \{T v \mid v \in [0, y]\} =
\\
y_0 + \rktp {T} {x} + \rktp {T} {y},
\end {multline*}
which implies that
\begin {equation}
\label {tztuv3}
\rktp {T} {x + y} \geqslant T z \geqslant y_0 + \rktp {T} {x} + \rktp {T} {y}.
\end {equation}
However, $\rktps {T}$ is linear by the assumptions, so
$$
\rktp {T} {x + y} = \rktp {T} {x} + \rktp {T} {y},
$$
which contradicts
\eqref {tztuv3} since $y_0 > 0$.  The proof of Proposition~\ref {lrdp} is complete.

\section {Dominating operators for superlinear maps}

\label {dofsm}

In this section we will establish a lemma in two versions that is essential for the main results of this paper,
and we will also verify some auxiliary propositions.
It is not difficult to see that (at least under some additional assumptions)
the validity of this lemma for Riesz-Kantorovich transforms in place of superlinear maps $p$
is more or less equivalent to validity of the main results in
Section~\ref {rkf} below concerning lattice and decomposition properties of a space of linear operators between two ordered linear spaces.

\begin {lemma}
\label {suprtli}
Suppose that $X$ is an ordered linear space such that
$x \in X_+$ is an internal point of $X_+$ and $Y$ is a Dedekind complete lattice.
Suppose also that $p : X_+ \to Y$ is a superlinear map and $y \in Y$ satisfies $y \geqslant p (x)$.
Then there exists a linear operator $M : X \to Y$ such that $M \geqslant p$ and $M x = y$.
\end {lemma}

Indeed, observe that if $x$ is an internal point of $X_+$ (equivalently, an order unit of $X$)
then the space $L = \{\alpha x \mid \alpha \in \mathbb R\}$ spanned by
$x$ satisfies $X_+ + L = X$. 
Let $M_0 : \{\lambda x \mid \lambda > 0 \} \to Y$ be a map defined by $M_0 (\lambda x) = \lambda y$ for all $\lambda > 0$.
Then $M_0 \geqslant p$ on the domain of $M_0$.
We see that $M_0$ satisfies the conditions of the Hahn-Banach theorem (see Section~\ref {ovl}),
and thus $M_0$ can be extended to a linear map $M$ satisfying the conclusion of Lemma~\ref {suprtli}.

\begin {comment}
This gives the following corollary to the Hahn-Banach theorem
which will be used later in Section~\ref {dofsm}.
\begin {theorem}
\label {onedimext}
Suppose that $X$ and $Y$ are ordered linear spaces,
$Y$ is Dedekind complete,
$p : X_+ \to Y$ is a superlinear map, $e \in X_+$ is an order unit and $y \in Y$ satifsfies $y \geqslant p (e)$.
Then there exists $T \in \linop {X} {Y}$ such that $T \geqslant p$ and $T e = y$.
\end {theorem}
\end {comment}

\begin {lemma}
\label {suprtlf0}
Suppose that $X$ is an ordered linear space equipped with a locally convex topology (not necessarily Hausdorff)
such that the cone $X_+$ is generating for $X$.
Suppose also that $p : X_+ \to \mathbb R$ is a continuous superlinear map and points $x \in X_+$, $y \in \mathbb R$ satisfy $y > p (x)$.
Then there exists a continous linear functional $M : X \to \mathbb R$ such that $M \geqslant p$ and $M (x) < y$.
\end {lemma}

This lemma is a little more involved since we cannot in general just invoke the Hahn-Banach theorem in the absence of interior points;
however, with a bit of topological reasoning we can instead make use of a separation theorem.
Let
$\Gamma = \{ (z, t) \mid z \in X_+, t \leqslant p (x) \}$ be the subgraph of $p$, that is the epigraph of $-p$ turned upside down,
and let $G = \bar \Gamma$ be its closure in $X \times \mathbb R$.
Observe that since $p$ is continuous, $\Gamma$ is closed in $X_+ \times \mathbb R$ equipped with the relative topology induced by $X \times \mathbb R$.
Indeed, for every $z \in X_+$ and $t \neq p (z)$ there exists a neighbourhood $V \subset X_+$ of $z$ such that
$|p (a) - t| < \frac 1 2 |p (z) - t|$ for all $a \in V$; therefore $V \times (t - \frac 1 2 |p (z) - t|, t + \frac 1 2 |p (z) - t|)$
is an open neighbourhood of
$(z, t)$ that does not intersect the graph $\Gamma$.
Closedness of $\Gamma$ in $X_+ \times \mathbb R$ implies that $G \cap (X_+ \times \mathbb R) = \Gamma$:
inclusion $G \cap (X_+ \times \mathbb R) \supset \Gamma$ is trivial, and nonemptiness of
$\left[G \cap (X_+ \times \mathbb R)\right] \setminus \Gamma$ would contradict relative closedness of $\Gamma$ in $X_+ \times \mathbb R$
since every point of $G \cap (X_+ \times \mathbb R)$ belongs to the closure of $\Gamma$ in $X_+ \times \mathbb R$.
In particular, $(x, y) \notin G$ since $(x, y) \notin \Gamma$.
Therefore point
$(x, y)$ belongs to $(X \times \mathbb R) \setminus G$, which is an open set in the locally convex linear topological space $X \times \mathbb R$,
and so there exists a convex open neighbourhood $U \ni (x, y)$ that
does not intersect $G$.
Superlinearity of $p$ implies that $\Gamma$ is a convex cone,
and therefore $G$ is also a convex cone as a closure of a convex cone in the linear topological space $X \times \mathbb R$.
Let $W = \bigcup_{\lambda > 0} \lambda U$ be the conic hull of $U$. $W$ is open because $U$ is open.
Since $G$ is a cone and $U$ does not intersect $G$, it is easy to see
that $W$ also does not intersect $G$.
By the Hahn-Banach separation theorem there exists a linear functional $f : X \times \mathbb R \to \mathbb R$ that separates $G$ and~$W$,
that is
$s = \sup \{f (b) \mid b \in G \} < f (c)$ for all $c \in W$.  Since $W$ is open and nonempty, $f$ is continuous.
Now it remains to verify that $D = \{(z, t) \mid f ((z, t)) = s\}$ is a graph of a functional $M$ that satisfies the conclusion of Lemma~\ref {suprtlf0}.
First, observe that $\{f ((z, t)) \mid t \in \mathbb R\} = \mathbb R$ for all $z \in X_+$:
since by linearity $f ((z, t)) = f ((x, t)) + f ((z - x, 0))$, it is sufficient to establish that
$f (x, \mathbb R) = \{f ((x, t)) \mid t \in \mathbb R\} = \mathbb R$,
and this follows easily from the fact that $f (x, \mathbb R)$ is a linear subspace of $\mathbb R$ that cannot be trivial since
$f ((x, p (x))) < f ((x, y))$.  So $D$ is indeed a graph of a functional
$M : X \to \mathbb R$, which is continuous since $f$ is continuous and therefore $D$ is closed.
Linearity of $f$ implies that $M$ is affine, and from the fact that $f$ separates two cones it follows that $f (0, 0) = s$
and $M (0) = 0$, i.~e. $M$ is linear.  Since the epigraph of $M$ contains an open set $U$, the epigraph of $M$ has nonempty interior
and thus $M$ is continuous.
Finally, $M \geqslant p$ because $f ((z, p (z))) \leqslant s$ for all $z \in X_+$, and 
$M (x) < y$ because $f ((x, M (x))) = s < f ((x, y))$.  The proof of Lemma~\ref {suprtlf0} is complete.

In order to apply Lemma~\ref {suprtlf0} to the Riesz-Kantorovich transform $\rkts {T}$ we need to make sure that $\rkts {T}$ is continuous
as long as functionals $T$ are continuous.
\begin {proposition}
\label {rktc}
Suppose that $X$ is an ordered linear topological space.
Then for any continuous linear functionals $T = \{T_1, \ldots, T_n\}$,
where $T_j : X \to \mathbb R$, $1 \leqslant j \leqslant n$, the Riesz-Kantorovich transform
$\rkts {T}$, assuming it is well-defined, is also continuous.
\end {proposition}

Indeed, suppose that $\rkt {T} {x} = y$ for some $x \in X_+$ and $y \in \mathbb R$ and we are given an open neighbourhood
$(y - \varepsilon, y + \varepsilon) \subset U \subset \mathbb R$ of $y$, $\varepsilon > 0$;
we need to find an open neighbourhood $V \subset X_+$ of $x$ satisfying $\rkt {T} {V} \subset U$.
Let $W = \left(-\varepsilon, \varepsilon\right) \subset \mathbb R$, so that $y + W \subset U$.
Observe that a map $S : X^n \to \mathbb R$ defined by $S ((x_1, \ldots, x_n)) = \sum_{j = 1}^n T_j (x_j)$ is linear and continuous in $x_j \in X$,
$1 \leqslant j \leqslant n$, as a composition of continuous maps.
Moreover, $S$ is uniformly continuous because of its linearity, i.~e. there exists an open neighbourhood
$V_1$ of $0$ such that
\begin {equation}
\label {unicont}
S ((x_1 + V_1, \ldots, x_n + V_1)) \subset S ((x_1, \ldots, x_n)) + W
\end {equation}
for all $x_j \in X$, $1 \leqslant j \leqslant n$.
Now let $V \subset X$ be an open balanced (i.~e. $V = -V$) neighbourhood of $0$ satisfying $V \subset V_1$.
Inclusion \eqref {unicont} implies that
\begin {equation}
\label {tzx}
\rkt {T} {z_1} \geqslant \rkt {T} {z_0} - \varepsilon
\end {equation}
for all $z_0 \in X_+$ and $z_1 \in \left( z_0 + \sum_{k = 1}^n V \right) \cap X_+$.
By setting $z_0 = x$ and $z_1 = x$ separately in \eqref {tzx} we arrive at the required continuity relation
$|\rkt {T} {z} - \rkt {T} {x}| \leqslant \varepsilon$ for all $z \in x + \sum_{k = 1}^n V$.  The proof of Proposition~\ref {rktc} is complete.

\begin {comment}
Let us give an immediate application of Proposition~\ref {rktc} which will be used later in the proof of Theorem~\ref {suprt}.
It is well known that a concave function is continuous on an open interval, but may fail to be continuous at the endpoints of a segment,
although it is always lower semicontinuous.
Is the Riesz-Kantorovich transform continuous when it is restricted on a segment?  We say that for a linear space $X$ a set
$\mathcal L \subset \linop {X} {\mathbb R}$ \emph {separates points of $X$} of $f (x) = 0$ for all $f \in \mathcal L$ implies that $x = 0$.
\begin {proposition}
\label {rksc}
Suppose that $X$ is an ordered linear topological space such that $\blinop {X} {\mathbb R}$ separates points of $X$
and $Y$ is a Dedekind complete ordered linear space.
Let $\ell \subset X$ be a line.
Suppose also that $T = \{T_1, \ldots, T_n\}$
are order bounded operators $T_j : \blinop {X} {Y}$, $1 \leqslant j \leqslant n$ and $S \in \blinop {X} {Y}$ satisfies $S \geqslant \rkts {T}$.
Then the set
$$
A = \{z \in \ell \cap X_+ \mid \rkt {T} {z} < S z\}
$$
is relatively open in $L = \ell \cap X_+$ with the topology induced from
$\mathbb R$, i.~e. for any $x \in A$ and $z \in L$ there exists some $\varepsilon > 0$
such that $(1 - \theta) x + \theta z \in A$ for all $0 < \theta < \varepsilon$.
\end {proposition}
First, let us equip $X$ with its order topology.
Observe that $\rkts {T}$ is superlinear, so its restriction on $L$ is a concave map.
It follows that the set $\left(\ell \cap X_+\right) \setminus A = \{z \in \ell \cap X_+ \mid \rkt {T} {z} = S (z)\}$
is convex, and therefore $A$ consists of at most two convex subsets of $\ell$; we need to show that none of these subsets is a point.
Suppose, on the contrary, that for some $x \in A$ and $z \in L$ we have
$z_\theta = (1 - \theta) x + \theta z \notin A$ for all $0 < \theta \leqslant 1$.
By Theorem~\ref {posfuncded} there exists a positive linear functional $g$ on $Y$ satisfying
$g \left(S (x) - \rkt {T} {x}\right) = 1$.  Let $S_1 = S \circ g$ and $R = \rkts {T} \circ g$,
so that $S_1 (x) = R (x) + 1$ but $S_1 (z_\theta) = R (z_\theta)$ for all $0 < \theta \leqslant 1$.
Proposition~\ref {monsuprel} implies that
$R \geqslant \rkts {T \circ g}$.
Thus $S_1 (z_\theta) = R (z_\theta)$ for all $0 < \theta \leqslant 1$ but $S_1 (x) = R (x) + 1 \geqslant \rkt {T \circ g} {x} + 1$;
in particular, $S_1 (x) \neq \rkt {T \circ g} {x}$.
Observe that $S_1$ is order bounded and therefore continuous, and likewise $T_j \circ g$ are continuous linear functionals for $1 \leqslant j \leqslant n$,
so $R$ is also continuous by Proposition~\ref {rktc}.  Since
$\blinop {X} {\mathbb R}$ separates points of $X$, the order topology on $X$ is Hausdorff, so its restriction on $L$ is also Hausdorff,
which means that every continuous function $\kappa : L \to \{0, 1\}$ is constant.
Setting $\kappa = S_1 - R$ we arrive at a contradiction which concludes the proof of Proposition~\ref {rksc}.

Proposition~\ref {rksc} tells us that, roughly speaking, the behavior of the Riesz-Kantorovich transform of linear operators
at the non-internal points of
the ordering cone is fairly regular and resembles continuity.  One might further wonder whether the Riesz-Kantorovich transform of linear operators
is in some sense differentiable at every point of the cone.
Unfortunately it might not.  To see this, let $X = \mathbb R^3$ and
let $X_+ = \{(x, y, z) \mid y^2 + (z - x)^2 \leqslant x^2\}$ be the ice-cream cone such that the $x$-axis is its extreme ray.
Let $f ((x, y, z)) = y$.  It is not difficult to see that the Riesz-Kantorovich transform $\rkts {f}$
has infinite derivatives in the direction $(0, 0, 1)$ at all points on the positive $x$-axis.
\end {comment}

We will also need the following simple proposition for the proof of Theorem~\ref {grkf} below.
Recall that an ordered linear space $Y$ is called \emph {Archimedean} if for any $y \in Y$ and $x \in Y_+$ inequalities $n y \leqslant x$ for all
$n \geqslant 1$ imply that $y \leqslant 0$.
\begin {proposition}
\label {isameorder}
Suppose that $X$ and $Y$ are ordered linear spaces and $\tint X_+ \neq 0$.  Suppose also that $Y$ is an Archimedean Dedekind complete lattice.
Then $X$ ordered by $\{0\} \cup \tint X_+$ defines the same order in $\linop {X} {Y}$ as
does $X$ ordered by $X_+$.
\end {proposition}
Indeed, suppose that $T \in \linop {X} {Y}$ and $T x \geqslant 0$ for all $x \in \tint X_+$ under the assumptions of Proposition~\ref {isameorder};
we need to verify that $T x \geqslant 0$ for all $x \in X_+$.  Fix any $x_0 \in \tint X_+$ and let
$x_\theta = (1 - \theta) x_0 + \theta x$ for $0 \leqslant \theta \leqslant 1$.
Since $x_\theta \in \tint X_+$ for $0 \leqslant \theta < 1$, by the assumptions
$T (x_\theta) \geqslant 0$ for $0 \leqslant \theta < 1$, and it suffices to verify that $T (x_1) \geqslant 0$.
Since $Y$ is a Dedekind complete lattice, there exists
$y_\alpha = \inf \{T (x_\theta) \mid \alpha \leqslant \theta < 1\}$ for all $0 \leqslant \alpha < 1$ and $y_\alpha \geqslant 0$.
Linearity of $T (x_\theta)$ in $\theta$ implies that
$$
T (x_\alpha) \wedge T (x_1) \leqslant T (x_\theta) \leqslant T (x_\alpha) \vee T (x_1)
$$
for all
$\alpha \leqslant \theta \leqslant 1$,
so 
\begin {equation}
\label {veeest}
0 \vee [T (x_\alpha) \wedge T (x_1)] \leqslant y_\alpha \leqslant T (x_\alpha) \vee T (x_1)
\end {equation}
for all $0 < \alpha < 1$.
Now observe that linearity of $T (x_\theta)$ in $\theta$ also implies that the set-valued map
$\alpha \mapsto \{y \in Y \mid y \geqslant T (x_\alpha), y \geqslant T (x_1)\}$ is affine in $0 \leqslant \alpha \leqslant 1$,
so $\alpha \mapsto T (x_\alpha) \vee T (x_1)$ and likewise $\alpha \mapsto T (x_\alpha) \wedge T (x_1)$ are affine in $0 \leqslant \alpha \leqslant 1$;
therefore $\varphi (\alpha) = T (x_\alpha) \vee T (x_1) - T (x_\alpha) \wedge T (x_1)$ is an affine positive map and $\varphi (1) = 0$.
In particular,
$$
T (x_\alpha) \vee T (x_1) = T (x_\alpha) \wedge T (x_1) + \varphi (\alpha) \leqslant T (x_1) + \varphi (0),
$$
so
the set of values of the right-hand part of \eqref {veeest} has an upper bound.
Let
$$
m_1 = \sup \{T (x_\alpha) \wedge T (x_1) \mid 0 \leqslant \alpha < 1\},
$$
$$
y_1 = \sup \{y_\alpha \mid 0 \leqslant \alpha < 1\}
$$
and
$$
M_1 = \inf \{T (x_\alpha) \vee T (x_1) \mid 0 \leqslant \alpha < 1\}.
$$
It is easy to see that $0 \leqslant m_1 \leqslant y_1 \leqslant M_1$;
moreover,
$$
M_1 - m_1 \leqslant \varphi (\alpha) = (1 - \alpha) \varphi (0)
$$
for all $0 \leqslant \alpha < 1$, so
the Archimedean property of $Y$ implies that $M_1 - m_1 \leqslant 0$, and thus $0 \leqslant m_1 = y_1 = M_1$.
Furthermore,
$$
m_1 \leqslant T (x_1) \leqslant M_1,
$$
which implies that $m_1 = M_1 = T (x_1)$, and therefore
$y_1 = T (x_1) \geqslant 0$.  The proof of Proposition~\ref {isameorder} is complete.

We mention that Poposition~\ref {isameorder} may fail without the assumption that $Y$ is Dedekind complete;
it is easy to construct a finite-dimensional example with some closed $X_+$ and open $Y_+$.  It is, however, not clear whether the assumption
that $Y$ is Archimedean is necessary for the conclusion of Proposition~\ref {isameorder}.

\section {The Riesz-Kantorovich formula}

\label {rkf}

In this section we will show how the results of Section~\ref {dofsm} above can be applied to the question of
whether the supremum of a collection of order bounded linear operators coincides with its Riesz-Kantorovich transform
and explore the relationship between lattice and decomposition properties of an ordered linear space $X$ and the same properties of a space
$\mathcal L \subset \blinop {X} {Y}$ of order bounded linear operators between $X$ and $Y$.

As it was explained in Section~\ref {ovl},
there is a long-standing open problem that can be rephrased as follows:
is it possible to characterize linearity of the Riesz-Kantorovich transform of a particular
finite collection $T$ of linear operators in terms of the existence of the supremum of $T$ in an ordered space of linear operators?
Now we are going to establish some partial answers to this question.
First, let us treat the case of linear continuous functionals.
\begin {theorem}
\label {suprtc}
Suppose that $X$ is an ordered locally convex linear topological space, the cone $X_+$ is generating for $X$
and $\clinop {X} {\mathbb R} \subset \blinop {X} {\mathbb R}$.
Then for any $T = \{T_1, \ldots, T_n\} \subset \clinop {X} {\mathbb R}$ such that $\sup T$ exists in $\mathcal L = \clinop {X} {\mathbb R}$ we have
$\sup T = \rkts {T}$.
\end {theorem}
Suppose that under the assumptions of Theorem~\ref {suprt}
$T \subset \clinop {X} {\mathbb R}$
is a finite collection of functionals such that the supremum of $T$ exists in $\clinop {X} {\mathbb R}$.
Let us denote this supremum by $S = \sup T$; we need to verify that $S = \rkts {T}$.
Observe that from Proposition~\ref {rktinf} it follows that $\rkts {T} \leqslant S$.
Now suppose that $\rkts {T} = S$ does not hold true, which means that there exists some $x \in X_+$
such that $\rkt {T} {x} < S (x)$.
By Proposition~\ref {rktc} the Riesz-Kantorovich trasform $\rkts {T}$ is continuous, and since it is superlinear
we can apply Lemma~\ref {suprtlf0} to points $x$ and $y = S (x)$ with $p = \rkts {T}$ to obtain a linear continuous functional $M \geqslant \rkts {T}$
satisfying $M (x) < S (x)$.
However, $S = \sup T$ implies that $S \leqslant M$ since $M \geqslant \rkts {T} \geqslant T_j$ for all $1 \leqslant j \leqslant n$, and therefore
$M (x) \geqslant \sup T (x) = S (x)$.  This contradiction concludes the proof of Theorem~\ref {suprtc}.

There is one immediate application of Theorem~\ref {suprtc}.
\begin {theorem}
\label {suprtf}
Suppose that $X$ is an ordered linear topological space, the cone $X_+$ is generating for $X$,
and $\mathcal L$ is a linear space of linear functionals satisfying $\rlinop {X} {\mathbb R} \subset \mathcal L \subset \blinop {X} {\mathbb R}$.
Then for any finite collection $T = \{T_1, \ldots, T_n\} \subset \mathcal L$ such that $\sup T$ exists in $\mathcal L$ we have
$\sup T = \rkts {T}$.
\end {theorem}
Indeed, observe that the supremum $S = \sup T$ in $\mathcal L$ is also a supremum of $T$ in $\blinop {X} {\mathbb R}$,
since any upper bound $U \in \blinop {X} {\mathbb R}$ of $T$ also belongs to $\mathcal L$ because
$U \geqslant T_1$, so $U - T_1 \geqslant 0$ and $U$ is a sum of $T_1 \in \mathcal L$ and a regular operator $U - T_1$.
Equipping $X$ with its order topology (see Section~\ref {ovl}) reduces
Theorem~\ref {suprtf} at once to a
direct application of Theorem~\ref {suprtc}.

Now we a going to treat the case of order bounded operators on a cone with an internal point.
\begin {theorem}
\label {suprt}
Suppose that $X$ and $Y$ are ordered linear spaces, the cone $X_+$ is generating for $X$,
$X_+$ has an internal point, $Y$ is a Dedekind complete lattice
and $\mathcal L \subset \linop {X} {Y}$
is a linear space
of linear operators acting from $X$ to $Y$ satisfying
$
\rlinop {X} {Y} \subset \mathcal L \subset \blinop {X} {Y}.
$
Then for any $T = \{T_1, \ldots, T_n\} \subset \mathcal L$ such that $\sup T$ exists in $\mathcal L$ we have
$\sup T (x) = \rkt {T} {x}$ for all $x \in \tint X_+$.
\end {theorem}
The proof essentially repeats the proof of Theorem~\ref {suprtc} above with only a couple of small twists.
Suppose that 
$T \subset \mathcal L$ is a finite collection of operators under the assumptions of Theorem~\ref {suprt}
such that the supremum of $T$ exists in $\mathcal L$.
We denote this supremum by $S = \sup T$ and it is set upon us to verify that $S = \rkts {T}$ on the set of all internal points of $X_+$.
From Proposition~\ref {rktinf} it follows that $\rkts {T} \leqslant S$.
Suppose that $\rkt {T} {x} = S (x)$ does not hold true for all $x \in \tint X_+$.
Then there exists some $x \in \tint X_+$
such that $\rkt {T} {x} < S (x)$.
Let $y = \frac 1 2 \left( \rkt {T} {x} + S (x) \right)$; then
$\rkt {T} {x} < y < S (x)$.
By 
Lemma~\ref {suprtli}
there exists
a linear map $M \in \linop {X} {Y}$ such that
$M \geqslant \rkts {T}$ and $M (x) = y < S (x)$.  Observe that $M \geqslant T_1$, so $M - T_1 \in \rlinop {X} {Y}$ and therefore
$M \in \rlinop {X} {Y} + \mathcal L = \mathcal L$.
However, $S = \sup T$ in $\mathcal L$ implies that $S \leqslant M$, so
$M (x) \geqslant S (x)$.  This contradiction concludes the proof of Theorem~\ref {suprt}.

Unfortunately, it does not seem to be clear whether $\rkts {T}$ is linear on all $X_+$ under the conditions of Theorem~\ref {suprt}.
Luckily, however, this omission in our understanding of the Riesz-Kantorovich transform does not get in the way of the applications
connecting the properties of $X$ and $\mathcal L$, since the operators in $\mathcal L$ are uniquely defined by their values on $\tint X_+$.

Let us now observe that in the same three cases treated in Theorems~\ref {suprtc}--\ref {suprt}
above we can replace the assumption that
$\sup S$ exists in $\mathcal L$ by the assumption that $\mathcal L$ has the Riesz Decomposition Property
and still arrive at the corresponding linearity of $\rkts {T}$.
We will only prove it for the two cases that are used in Theorem~\ref {grkf} below.
\begin {proposition}
\label {suprtrdpf}
Suppose that $X$ is an ordered linear space, the cone $X_+$ is generating for $X$,
and $\mathcal L \subset \linop {X} {\mathbb R}$
is a linear space
of linear functionals on $X$ satisfying
$
\rlinop {X} {\mathbb R} \subset \mathcal L \subset \blinop {X} {\mathbb R}.
$
Suppose also that $\mathcal L$ has the Riesz Decomposition Property.
Then $\rkts {T} \in \mathcal L$ for any $T = \{T_1, \ldots, T_n\} \subset \mathcal L$.
\end {proposition}

\begin {proposition}
\label {suprtrdp}
Suppose that $X$ and $Y$ are ordered linear spaces, the cone $X_+$ is generating for $X$,
$X_+$ has an internal point, $Y$ is a Dedekind complete lattice
and $\mathcal L \subset \linop {X} {Y}$
is a linear space
of linear operators acting from $X$ to $Y$ satisfying
$
\rlinop {X} {Y} \subset \mathcal L \subset \blinop {X} {Y}.
$
Suppose also that $\mathcal L$ has the Riesz Decomposition Property.
Then $\rkts {T}$ is linear on $\tint X_+$ for any $T = \{T_1, \ldots, T_n\} \subset \mathcal L$,
and the unique extension of $\rkts {T}$ restricted to $\tint X_+$ to a linear operator from $\linop {X} {Y}$ belongs to~$\mathcal L$.
\end {proposition}

We will now give a proof for Proposition~\ref {suprtrdp} that also works as a proof of Proposition~\ref {suprtrdpf} with suitable apparent
minor modifications.
Indeed, under the assumptions of Proposition~\ref {suprtrdp} the Riesz-Kantorovich transform $\rkts {T}$ is well-defined and superlinear.
Suppose that it is not linear on $\tint X_+$; this means that
$\rkt {T} {x + y} > \rkt {T} {x} + \rkt {T} {y}$
for some $x, y \in \tint X_+$.
Let $z = \rkt {T} {x + y} - \left[\rkt {T} {x} + \rkt {T} {y} \right] > 0$.
By Lemma~\ref {suprtli} we can find some $S_0, S_1 \in \linop {X} {Y}$ such that
$S_0, S_1 \geqslant \rkts {T}$ and
\begin {equation}
\label {s0s1}
S_0 (x) = \rkt {T} {x} + \frac 1 4 z, \quad S_1 (y) = \rkt {T} {y} + \frac 1 4 z.
\end {equation}
By the same reasoning as in the proof of Theorem~\ref {suprt} we have $S_0, S_1 \in \mathcal L$.
Since by the assumption $\mathcal L$ has the Riesz Decomposition Property, it also has the interpolation property (see Section~\ref {ovl}).
Observe that $\{S_0, S_1\} \geqslant \rkts {T} \geqslant T$, so $\{S_0, S_1\} \geqslant T$, and by the interpolation property of $\mathcal L$
there exists some $S \in \mathcal L$
such that $\{S_0, S_1\} \geqslant S \geqslant T$.  Therefore $S \geqslant \rkts {T}$ by Proposition~\ref {rktinf}.
It follows that
\begin {multline}
\label {sestrk}
S (x + y) - [S (x) + S (y)] \geqslant \rkt {T} {x + y} -[S_0 (x) + S_1 (y)] =
\\
\left(\rkt {T} {x + y} - \left[\rkt {T} {x} + \rkt {T} {y} \right]\right) - \frac 1 2 z = z - \frac 1 2 z = \frac 1 2 z > 0,
\end {multline}
so $S$ cannot be linear.  This contradiction concludes the proof of Proposition~\ref {suprtrdp}.
In order to prove Proposition~\ref {suprtrdpf} we modify this proof as follows.
Surely we need to take any $x, y \in X_+$ and not just $x, y \in \tint X_+$.  As in the proof of Theorem~\ref {suprtf}
we equip $X$ with its order topology.
Then Lemma~\ref {suprtrdpf} instead of Lemma~\ref {suprtli} gives $S_0$ and $S_1$ satisfying inequalities $<$ in \eqref {s0s1}
rather than equalities, and thus we get inequality $>$ instead of the first equality in \eqref {sestrk}.
The rest of the proof does not change.

Now we can state an exact form of the Riesz-Kantorovich theorem that integrates the two general settings considered above.
\begin {theorem}
\label {grkf}
Suppose that $X$ and $Y$ are ordered linear spaces, $X_+$ is a generating cone in $X$,
$Y$ is a Dedekind complete lattice and $\mathcal L$ is a linear space of operators from $X$ to $Y$
satisfying $\rlinop {X} {Y} \subset \mathcal L \subset \blinop {X} {Y}$.
Suppose also that either the following assumptions is satisfied.
\begin {enumerate}
\item
$Y = \mathbb R$; we set $X_r = X_+$ in this case.
\item
$\tint X_+ \neq \emptyset$ and $Y$ is Archimedean; we set $X_r = \{0\} \cup \tint X_+$ in this case.
\end {enumerate}
Then the following conditions are equivalent.
\begin {enumerate}
\item
$\mathcal L$ is a lattice.
\item
$\mathcal L$ has the Riesz Decomposition Property.
\item
The Riesz-Kantorovich transform $\rkts {T}$ of any finite collection~$T$ of operators from $\mathcal L$ is linear on $X_r$.
\item
For any operator $T \in \mathcal L$ the Riesz-Kantorovich transform $\rktps {T}$ is linear on $X_r$.
\item
$X$ ordered by the cone $X_r$ has the $\mathcal L$-Riesz Decomposition Property.
\end {enumerate}
Moreover, if any (and therefore all) of these conditions are satisfied then $\mathcal L = \rlinop {X} {Y}$
and $\sup T = \rkts {T}$ in $\mathcal L$ on $X_r$ for any finite collection of operators $T \subset \mathcal L$.
\end {theorem}
It is well known that implication $1 \Rightarrow 2$ is satisfied for any ordered linear space $\mathcal L$.
Implication $2 \Rightarrow 3$ follows from either Proposition~\ref {suprtrdpf} or Proposition~\ref {suprtrdp}.
Implication $3 \Rightarrow 4$ is trivial.
Observe that the cone $X_r$ is generating for $X$,
and by Proposition~\ref {isameorder} space $X$ ordered by $X_r$ generates the same order in
$\linop {X} {Y}$ as does $X$ ordered by $X_+$.
Proposition~\ref {lrdp} proves the implication $4 \Rightarrow 5$.
Suppose now that Condition~5 of Theorem~\ref {grkf} is staisfied.
Then Proposition~\ref {lrdptorkf} implies that for any $T \in \mathcal L$ we have $\rktps {T} \in \mathcal L$,
and from Proposition~\ref {rktinf} it follows that $\rktps {T} = T^+ = T \vee 0$ in $\mathcal L$,
so $\mathcal L$ is a lattice (since then for any $S, T \in \mathcal L$ expression $(T - S)^+ + S$ provides the supremum $S \vee T$ of $\{S, T\}$).
We have thus verified implication $5 \Rightarrow 1$.
Finally, if Condition~1 is satisfied
then for any $T \in \mathcal L$ we have a decomposition $T = T \vee 0 - (-T) \vee 0$ into a difference of two positive linear maps that belong to
$\mathcal L$, so $\mathcal L \subset \rlinop {X} {Y}$, and the converse inclusion was assumed.
Thus $\mathcal L = \rlinop {X} {Y}$.
Of course, Condition~1 means that $\sup T$ exists in $\mathcal L$ for a finite $T \subset \mathcal L$,
so $\sup T = \rkts {T}$ follows from either Theorem~\ref {suprtf} or Theorem~\ref {suprt}.
The proof of Theorem~\ref {grkf} is complete.


\begin {comment}
Let us state a question related to the Riesz-Kantorovich formula that appears to  be still mysterious in the following more general form.
\begin {problem}
Suppose that $X$ and $Y$ are ordered linear spaces, $\blinop {X} {\mathbb R}$ separates points of $X$ and the cone $X_+$ is generating for $X$.
Assume that $T, S \in \blinop {X} {Y}$ and
let $P : X_+ \to Y$ be a superlinear map satisfying $T \leqslant P \leqslant S$;
suppose also that restrictions of $P$ on lines in $X$ are continuous in the sense that it satisfies the conclusion of Proposition~\ref {rksc}.
If
$
S = \sup \left\{M \in \blinop {X} {Y} \mid M \geqslant P\right\},
$
does it follow that $P = S$?
\end {problem}
In the case that $X$ is finite dimensional, it is fairly easy to give a positive answer to this problem
by repeating the proof of Theorem~\ref {suprt}.  However, this problem appears to be interesting even
in the case
$$
X = c_{00} = \left\{\{x_j\}_{j \in \mathbb N} \mid \text {$x_j = 0$ except for a finite number of indices $j$}\right\}
$$
of finite sequences with the standard ordering $X_+ = \{x_j \geqslant 0\}$.
\end {comment}



\bibliographystyle {plain}

\bibliography {bmora}

\end {document}

%% file: lwen.tex
\usepackage {hyperref}

\usepackage [utf8] {inputenc}

\usepackage{mathtext}

\usepackage {version}

\usepackage {amsmath}
\usepackage {amssymb}
\usepackage {amsthm}

\usepackage {graphicx}

\usepackage {appendix}

\usepackage{mathtext}
\usepackage{amsmath}
\usepackage{amsfonts}
\usepackage{amssymb}
\usepackage{mathrsfs}
\usepackage{amsthm}

\usepackage {needspace}

\excludeversion {svntags}

\newtheorem {theorem} {Theorem}
\newtheorem {problem} [theorem] {Problem}

\newtheorem* {theoremhahnbanach} {Hahn-Banach Theorem}
\newtheorem* {theoremrieszkantorovich} {Riesz-Kantorovich Theorem}
\newtheorem {lemma} [theorem] {Lemma}
\newtheorem {proposition} [theorem] {Proposition}
\newtheorem {definition} [theorem] {Definition}

\newcommand {\lclassg} [1] {\ensuremath{\mathrm L_{#1}}}

\DeclareMathOperator* {\tint} {int}

\excludeversion {comment}

\newcommand {\linop} [2] {\ensuremath{\mathfrak L (#1, #2)}}
\newcommand {\suplinop} [2] {\ensuremath{\mathfrak S (#1, #2)}}
\newcommand {\blinop} [2] {\ensuremath{\mathfrak L^b (#1, #2)}}
\newcommand {\orderdual} [1] {\ensuremath{{#1}^{\sim}}}
\newcommand {\rlinop} [2] {\ensuremath{\mathfrak L^r (#1, #2)}}
\newcommand {\clinop} [2] {\ensuremath{\mathfrak L^c (#1, #2)}}

\newcommand {\rkt} [2] {\ensuremath{\mathcal R (#1; #2)}}
\newcommand {\rkts} [1] {\ensuremath{\mathcal R (#1)}}
\newcommand {\rktp} [2] {\ensuremath{\mathcal R^+ (#1; #2)}}
\newcommand {\rktps} [1] {\ensuremath{\mathcal R^+ (#1)}}